\documentclass[11pt,reqno]{amsart}
\usepackage{amsmath,amsfonts,amssymb}

\usepackage[utf8]{inputenc}
\usepackage[T1]{fontenc}

\newlength{\aux}
\newtheorem{lemma}{Lemma}[section]
\newtheorem{property}{Property}[section]
\newtheorem{hyp}{Hypothesis}[section]

\newtheorem{theo}[lemma]{Theorem}

\newtheorem{proposition}[lemma]{Proposition}
\numberwithin{equation}{section}

\newcommand{\RR}{{\mathbb R}}
\newcommand{\NN}{{\mathbb N}}
\newcommand{\CC}{{\mathbb C}}
\newcommand{\PP}{{\mathbb P}}

\newcommand{\Id}{{\mbox{Id }}}
\newcommand{\Hess}{{\mbox{Hess }}}
\newcommand{\osc}{{\mbox{osc }}}

\newcommand{\tL}{{\widetilde L}}

\newcommand{\p}{{\mathcal P}}
\newcommand{\G}{{\mathcal G}}
\newcommand{\E}{{\mathbb E}}
\newcommand{\ppq}{\leqslant}
\newcommand{\pgq}{\geqslant}
\newcommand{\tr}{tr}
\newcommand{\mun}{{\widehat \mu_N}}

\newcommand{\nun}{{\widehat \nu_N}}
\newcommand{\rhon}{{\widehat \rho_N}}
\newcommand{\eps}{\varepsilon}
\newcommand{\un}{1}
\newcommand{\Lip}{{Lip}}
\newcommand{\hf}{{\widetilde f}}
\newcommand{\Ima}{\Im\mbox{m}}
\newcommand{\Rea}{\Re\mbox{e}}

\newcommand{\PNVb}{\PP^N_{V,\beta}}
\newcommand{\PNmVb}{\PP^{N-1}_{V,\beta}}

\newcommand{\cVb}{c_{V,\beta}}

\newcommand{\DVb}{D_{V,\beta}}

\newcommand{\PNV}{\PP^N_V}
\newcommand{\ZNV}{Z^N_V}
\newcommand{\ZNVb}{Z^N_{V,\beta}}

\newcommand{\ZNmVb}{Z^{N-1}_{V,\beta}}
\newcommand{\muVb}{\mu_{\frac{2V}{\beta}}}
\newcommand{\rVb}{\rho_{\frac{2V}{\beta}}}

\newcommand{\AVb}{A_{V,\beta}}

\newcommand{\SV}{\Sigma_V}

\newcommand{\StV}{{\widetilde \Sigma_{V}}}
\newcommand{\StVb}{{\widetilde \Sigma_{\frac{2}{\beta}V}}}

\def\Xint#1{\mathchoice
{\XXint\displaystyle\textstyle{#1}}
{\XXint\textstyle\scriptstyle{#1}}
{\XXint\scriptstyle\scriptscriptstyle{#1}}
{\XXint\scriptscriptstyle\scriptscriptstyle{#1}}
\!\int}
\def\XXint#1#2#3{{\setbox0=\hbox{$#1{#2#3}{\int}$ }
\vcenter{\hbox{$#2#3$ }}\kern-.5\wd0}}

\def\dashint{\Xint-}

\newenvironment{dem}{\textbf{Proof.}\par}
{\begin{flushright}$\Box$\end{flushright}}

\title[Free transport-entropy inequalities and application]{Free transport-entropy inequalities for non-convex potentials and application to concentration for random matrices}

\author[M. Maïda, \'E. Maurel-Segala]{Mylène Maïda$^\star$, \'Edouard Maurel-Segala$^\sharp$}
\thanks{Université Paris Sud 11
Laboratoire de Mathématiques, Bat. 425
91405 Orsay Cedex, France \\
$^\star$  E-mail: Mylene.Maida@math.u-psud.fr,\\
$^\sharp$ E-mail: Edouard.Maurel-Segala@math.u-psud.fr.\\
This work was   supported by the \emph{Agence Nationale de la
Recherche} grant ANR-08-BLAN-0311-03}

\date{\today}
\begin{document}

\maketitle
\begin{abstract}
Talagrand's inequalities provide a link between two fundamentals concepts of probability:  transportation
and  entropy.
The study of the  counterpart of these inequalities in the context of free probability has been initiated by Biane and 
Voiculescu and later extended by Hiai, Petz and Ueda for convex potentials.
In this work, 
we prove a free analogue of a result of Bobkov and Götze in the classical setting, thus providing free transport-entropy
inequalities for a very natural class of measures appearing in random matrix theory.
These inequalities are weaker than the ones of Hiai, Petz and Ueda but
 still hold beyond the convex case.
We then use this result to get a concentration estimate for $\beta$-ensembles under  mild assumptions on the potential.

\end{abstract}

\newpage

\section{Introduction}
\subsection{Classical transport-entropy inequalities}
In transportation theory, 
an important achievement was the proof by Talagrand in \cite{Tal} of the fact that the
standard Gaussian measure $\gamma$ in $\RR^n$ satisfies the transport-entropy inequality $T_2(2)$ (named after Talagrand).
We say that a probability measure $\mu$ on $\RR^n$ satisfies the inequality $T_p(C)$ for some $C>0$ if, 
for any probability measure $\nu$ on $\RR^n,$ 
$$ W_p^2(\nu, \mu) \ppq C H(\nu | \mu),$$ 
where
\begin{itemize}
\item
$W_p(\mu,\nu)$ is the Wassertein distance of order $p$ with respect to the 
Euclidean distance on $\RR^n$ between the two probability measures $\mu$ and $\nu,$ that is  
$$W_p(\mu, \nu) := \inf\left\{ \int_{(\RR^n)^2} |x-y|^p d\pi(x,y); \pi_0 =\mu, \pi_1=\nu\right\}^{\frac{1}{p}},$$
with $\pi_0$ and $\pi_1$ respectively the first and second marginals of $\pi$.
\item
 $H( \cdot | \mu)$ is the (classical) relative entropy  with respect to $\mu,$ that is
$$ H(\nu | \mu) = \int \ln\left(\frac{d\nu}{d\mu}\right) d\nu,$$ if $\nu$ is absolutely continuous with respect to $\mu$
and $+\infty$ otherwise.
If $\mu$ is the Lebesgue measure on $\RR^n,$   $H(\cdot |\mu)$ is just called the classical entropy and denoted by $H(\cdot).$
\end{itemize}

These inequalities give very important informations on measures that satisfy them since they are related to concentration properties                                                                                                                                        
and allow to deduce precise deviation estimates starting from a large deviation principle (see the work of Gozlan and Léonard \cite{GL} for a 
discussion on these topics as well 
as an excellent review of the advances during the past decade in this field).

\par\medskip

After the result of Talagrand, a lot of attention was devoted to prove similar inequalities beyond the Gaussian case; we will review only a 
few of them.
It was proved by Otto and Villani in \cite{OV} that any probability measure $\mu$ satisfying a log-Sobolev inequality with constant $C$ also 
satisfies the inequality $T_2(C)$. In particular, let $\mu$ be a probability measure
of the form $d\mu(x)= e^{-V(x)}dx,$ for some potential $V$ satisfying $\Hess V\pgq\kappa \Id>0.$
In that case,  for any probability measures $\nu$ on $\RR^n,$ 
\begin{equation}
 \label{blo} W_2^2(\nu, \mu) \ppq \frac 2 \kappa  H(\nu | \mu) 
\end{equation}
as if $\Hess V\pgq\kappa \Id>0,$ then the measure $d\mu(x)= e^{-V(x)}dx$ satisfies a log-Sobolev inequality with constant $2/\kappa$.

There has been some attempts (e.g. \cite{CG,CGW}) to generalise this results to potentials $V$ that are no longer strictly convex but 
the criteria that have been obtained are quite difficult to handle.

Furthermore, it seems that there is little room for improvements 
of the result of Otto and Villani since  the inequality $T_2$ implies Poincaré inequality for $\mu$.
Thus it is impossible to hope  for a measure $\mu$ that does not have a connected support 
to satisfy an inequality $T_2$ since such measures does not satisfy Poincaré inequality.
For such measures, one can be interested in the inequality $T_1$.
Note that 
this inequality $T_1(C)$ is weaker than $T_2(C)$ by a direct application of Cauchy-Schwarz inequality.
The benefit 
is that the criteria for $T_1$ are much easier to handle. In particular, in \cite{BG1},
Bobkov and Götze proved that a probability measure $\mu$ satisfies $T_1(2C)$ if and only if,
$$\int e^{f(x)} d\mu(x) \ppq e^{C\frac{\|f\|_\Lip^2}{2}}$$
for all Lipschitz function $f$ such that $\int f d\mu=0$ (with $\|f\|_\Lip$ denoting the Lipschitz constant of $f$).
Later, Djellout, Guillin and Wu proved in \cite{DGW} that this condition was equivalent to the quite easy to handle condition that there 
exists $\alpha>0$ and $x_0$
such that $$\int \exp(\alpha d(x,x_0)^2)d\mu(x)<+\infty.$$
One can see on this latter expression that compactly supported measures automatically satisfy a $T_1$ inequality. Besides, if $\mu$ is a measure 
of density $\exp(-V(x))$ with $V(x)\sim |x|^d$ for large $x$, then $\mu$ satisfies $T_1$ if and only if $d\pgq 2$ (note the similar condition 
appearing in Hypothesis \ref{hyp} below).

\subsection{Free transport-entropy inequalities}

We review hereafter some results in the literature that are the analogues in the free probability context of the inequality $T_2$
previously discussed. We assume that the reader has some minimal background in free probability, that can be found for example in \cite{AGZ}.

In the free probability context, the semi-circle law, also called Wigner law, given by 
$d\sigma(x)= \frac{1}{2\pi} \sqrt{4-x^2} \mathbf{1}_{[-2,2]}(x)dx$ can for many reasons be seen as the free analogue of the standard Gaussian 
distribution.
Therefore it is  natural to ask whether the semi-circle law satisfies a free analogue of the transport-entropy inequality $T_2,$
with the entropy  replaced by the free entropy defined by Voiculescu (see \cite{V} for a quick review).
A positive answer to this question was given by Biane and Voiculescu in \cite{BV}~: they showed that for any compactly supported probability 
measure $\nu,$
$$ W_2^2(\nu, \sigma) \ppq 2 \Sigma(\nu),$$
where $\Sigma$ is the free entropy with respect to $\sigma$ (called free entropy adapted to the free Ornstein-Uhlenbeck process 
in \cite{BV}).

The free entropy was introduced in whole generality (even for multivariate tracial states) by Voiculescu, it is a profound and quite complicated object 
but luckily in the one dimensional setting, one can give the following explicit expression for the free entropy with respect to $\sigma$~:
\begin{equation} \label{def:sigma}
 \Sigma(\nu) = \frac{1}{2} \int x^2 d\nu(x) - \iint \ln|x-y| d\nu(x)d\nu(y) - \frac{3}{4}.
\end{equation}

As we said that the semi-circle law $\sigma$ is  the analogue of the Gaussian law, one can now wonder what are the free analogues
 $\mu_V$ of the classical measures of the form $ e^{-V(x)}dx.$ To define those probability measures $\mu_V,$ we need to look at the probability 
measures defined on
the space of $N$ by $N$ Hermitian matrices by:
$$d\mu^N_V(X)\propto\exp(-N\tr V(X))d^NX$$
where $d^NX$ is the Lebesgue measure on the space of Hermitian matrices. 
In the sequel, we will assume that the potential $V$ satisfies  
\begin{hyp}\label{hyp}
 $V$ is continuous and $\liminf_{|x| \rightarrow \infty} \frac{V(x)}{x^2} >0$.
\end{hyp}
It ensures for example the existence of a normalising constant such that $\mu^N_V$ becomes a probability measure.
Note that this hypothesis is a little more restrictive than the usual growth requirement for this model but seems necessary for our result.
If the matrix $X^N$ is distributed according to the law $\mu^N_V$ then
the joint law of the eigenvalues of $X^N$ is the following~:
$$\PNV(dx_1,\dots,dx_N)= \prod_{i<j}|x_i -x_j|^2\exp\left(-N\sum_{i=1}^NV(x_i)\right)\frac{\prod_{i=1}^Ndx_i}{\ZNV},$$
with $\ZNV$ a normalising constant.
This can be seen as the density of a Coulomb gas, that is  
$N$ particles in the potential $NV$ with a repulsive electrostatic interaction.
Under the law $\PNV$, the particles $x_1,\dots,x_N$ tend to be near the minima of $V$ but due to the Vandermonde determinant they can not be too 
close from each other.
The study of how these two effects reach an equilibrium is a difficult, yet well studied one. We recall hereafter a few facts about their behaviour.
First, if we introduce the empirical measure $\mun:=\frac{1}{N}\sum_{i=1}^N\delta_{x_i}$, the density of $\PNV$ can be written as
$$\PNV(dx_1,\dots,dx_N)= \exp(-N^2{\widetilde J_V}(\mun))\frac{\prod_{i=1}^Ndx_i}{\ZNV}$$
with, for any probability measure $\mu,$
$$\widetilde J_V(\mu) = \int V(x)d\mu(x)-\iint_{x\neq y} \ln|x-y| d\mu(x) d\mu(y).$$

One can expect that in the large $N$ limit, the eigenvalues   should organise according to the minimiser of this functional.
We  recall hereafter a result of the classical theory of logarithmic potentials which will define the family of measures $\mu_V$ 
which are the analogues in the free probability setting to the probability measures of the form $ e^{-V(x)}dx.$
This result is Theorem 1.3 in Chapter 1 of \cite{ST} simplified by the use of Theorem 4.8 in the same chapter which implies the continuity of the 
logarithmic potential.
The books \cite{AGZ} and \cite{D} also give presentation of similar results, in a perspective closer to random matrix theory
but later on we will need some more involved results of the book of Saff and Totik so we try not to drift too 
much away from their notations.
Let us denote, for $X$ a Polish space, by  $\p(X)$ the set of probability measures on $X.
$
\begin{theo}[Equilibrium measure of a potential]\label{theo:st}
Let $V$ be a function satisfying Hypothesis \ref{hyp}.
Define for $\mu$ in $\p(\RR)$,
$$J_V(\mu) = \int_\RR V(x)d\mu(x)- \iint_{\RR^2} \ln|x-y| d\mu(x) d\mu(y)$$
with the convention $J_V(\mu)=+\infty$ as soon as $\int Vd\mu=+\infty$.
Then $c_V = \inf_{\nu \in \p(\RR)} J_V(\nu)$ is a finite constant and the minimum of $J_V$ is achieved at a unique probability measure $\mu_V$ called 
equilibrium measure which has a compact support.
Besides, if we define the logarithmic potential of $\mu_V$ as
$$ U_{\mu_V}(x) = -   \int \ln|x-y| d\mu_V(y),$$
for all $x \in \CC$ then $U_{\mu_V}$ is finite and continuous on $\CC$ and $\mu_V$ is the unique probability measure on $\RR$ for which there 
exists a constant $C_V$ such that:
$$
\begin{array}{rl}
-2U_{\mu_V}(x)+C_V\ppq V(x)&\mbox{for all $x$ in $\CC$.}\\
 -2U_{\mu_V}(x)+C_V=V(x)&\mbox{for all $x$ in the support of $\mu_V$}\\
\end{array}
$$
$C_V$ is related to $c_V$ by the formula $C_V = 2c_V-\int V(x) d\mu_V(x)$.
\end{theo}

This allows to define the free entropy relative to the potential $V$ as follows~: for any $\mu \in \p(\RR),$
$$\Sigma_V(\mu)=J_V(\mu)-c_V = J_V(\mu)-J_V(\mu_V).$$
This quantity is always positive and  vanishes only at 
$\mu_V$. One can check that the functional $\Sigma$ introduced in \eqref{def:sigma}
coincides with $\Sigma_{x^2/2}.$
\par\medskip

Let us make a few remarks on the functional $\Sigma_V.$
First, Theorem \ref{theo:st} studies the optimum for the functional $J_V$ but not how it is related to a typical distribution of $x_i$'s 
under the law $\PNV$.
This is the goal of the work of Ben Arous and Guionnet \cite{BaG} (see also the book \cite{AGZ} for a slightly different point of view), 
from which we want to recall the following result, that will be useful in the sequel.

\begin{theo}[Large deviations for the empirical measure]\label{theo:gd}
Let $V$ be a function satisfying Hypothesis \ref{hyp}.
Under the law $\mu^N_V,$ 
the sequence of random measures $\mun$ satisfies a large deviation principle in the speed $N^2$ with good rate function $\Sigma_V$.
\end{theo}
We refer the reader not familiar with the theory of large deviations to \cite{DZ}.
By comparison to Sanov theorem where the classical relative entropy appears as a good rate function, this result can be seen as a 
justification of the name ''free relative entropy''
for $\Sigma_V$.

Another reason is that $\Sigma_V$ appears as a limit of classical relative entropy. Indeed, under some additional assumptions on $V$
and $W,$ we have
$$\lim_N \frac{1}{N^2}H(\mu^N_W|\mu^N_V)=\Sigma_V(\mu_W).$$
A precise statement and a proof of this convergence will be given within the proof of Proposition \ref{theosup}
where it is needed.

We can now state a generalisation of the result of Biane and Voiculescu, which  can be seen as a free analogue 
of the classical result \eqref{blo}. It was first proved by Hiai, Petz and Ueda in \cite{HPU} using 
random matrix approximations   and   classical inequalities.
 Let $V$ be a strictly convex function with $V''(x)\pgq\kappa>0$ on $\RR$.
Then, for any probability measure $\nu$
$$ W_2^2(\nu, \mu_V) \ppq \frac{2}{\kappa} \Sigma_V(\mu).$$
The same result was later proved in a very direct way by Ledoux and Popescu \cite{LP}.

Finally, let us finish this quick review by mentioning two interesting directions that could extend
these works.
First,
in view of this result and the one by Otto and Villani, a natural question is to ask whether 
a free analogue of the log-Sobolev inequality (see the work of Biane \cite{B} for the construction of such an object) is sufficient to 
obtain a free transport inequality.
While the methods of \cite{BV} have some similarities with the ones of \cite{OV} this remains an open problem.

Another natural extension of these results would be to look at the multivariable case.
As pointed out above, in several variables, the free entropy is a much more difficult object to handle and the theory of 
non-commutative transport is still at its beginning. 
The recent paper of Guionnet and Shlyakhtenko \cite{GS} gives some basis and highlights many pitfalls of this  theory. 
Still, the Wasserstein distance is still well defined and in some cases such as a $n$-uple of semi-circular variables one can define 
a notion of free relative entropy.
In \cite{BD}, Biane and Dabrowski prove 
 a version of the free Talagrand inequality for a $n$-uple of
semi-circular variables.

\subsection{Statement of the free $T_1$ inequality}
The problem we want to address in this work is to prove
a free analogue of the result of Bobkov and Götze, thus providing free transport-entropy inequality for measures $\mu_V$ beyond the 
case of convex potentials which was treated in the work of Hiai, Petz and Ueda. As pointed out above, even in the 
classical context, there is no reason for measures coming non-convex potentials to satisfy $T_2$. 
Thus we will prove an analogue to the inequality $T_1.$ Our main result can be stated as follows
\begin{theo}[Free $T_1$ inequality] \label{main}
 Let $V$ be a function satisfying Hypothesis  \ref{hyp}. Then there exists a constant $B_V$ 
such that, for any probability measure $\nu$ on $\RR$,
$$ W_1^2(\nu,\mu_V) \ppq B_V \Sigma_V(\nu).$$
\end{theo}

Let us make a quick remark on the role on Hypothesis \ref{hyp}.
It is not hard to check that the result is trivially false if $V$ is negligible  with respect to $x^2$.
Indeed, if $\nu_n$ is the uniform law on $[n;n+1],$ as $\mu_V$ is compactly supported, $W_1(\nu_n,\mu_V)^2$ is equivalent to $n^2$ 
but $\Sigma_V(\nu_n)$ grows like $\nu_n(V)$ which would be less than quadratic.

A natural strategy to try to prove this theorem, following the idea \cite{HPU}, is to look at a
finite dimensional approximation by matrix models.
The issue with this approach is that while for the classical $T_2$ inequality the constant in front of the entropy is explicitly related 
to the potential and behave nicely when the dimension
increases, this is no longer the case for $T_1.$ In \cite{BVi}, Bolley and Villani managed
to explicitly link the constant to the potential but when applied in this case the constant deteriorates very quickly with the dimension. 
Thus we will need some new tools to get our results.
The main ingredients that we will use to adapt the proof of Bobkov and Götze
is potential theory.\par\medskip

Since Theorem \ref{main} is only stated for measures of the form $\mu_V$, one may have the false impression 
that it is restricted to this particular case.
In fact it is relevant for  a quite large class of measures.
A difficulty is that one may want to think of the functional $\Sigma_V$ as the entropy relative to the measure $\mu_V$ but we must be 
careful since different $V$'s can lead
to the same equilibrium measure while defining
different notions of this relative entropy.

The first step is to get rid of this dependence on the potential.
Let $\mu$ be a probability measure with a compact support $S_\mu$ in $\RR$ such that its logarithmic potential
$ U_\mu(x) = -   \int \ln|x-y| d\mu(y)$ exists and is continuous on $\CC$.
Then the potential $V(x)=-2U_\mu(x)+(d(x,S_\mu))^2$  satisfies Hypothesis \ref{hyp} and
using Theorem \ref{theo:st} it is easy to see that $\mu_V=\mu$.

Now if we look at $\nu$ a probability measure on $S_\mu$:
\begin{align*}
\Sigma_V(\nu)&=\int V d\nu - \iint_{\RR^2} \ln|x-y| d\nu(x) d\nu(y)-c_V\\
&=2\iint_{\RR^2}  \ln|x-y| d\mu(x) d\nu(y) - \iint_{\RR^2} \ln|x-y| d\nu(x) d\nu(y)-c_V+C_V\\
&=-\iint_{\RR^2} \ln|x-y| d(\nu-\mu)(x) d(\nu-\mu)(y)
\end{align*}
where we used the Theorem \ref{theo:st} on the second line and it is easy to check that there is no constant in the last line since the expression
must be $0$ for $\nu=\mu$.

This allows to define a relative free entropy which does not depend on a potential but only on a measure:
$$\Sigma_\mu(\nu)=-\iint_{\RR^2} \ln|x-y| d(\nu-\mu)(x) d(\nu-\mu)(y)$$
if $\nu$ has a support included in $S_\mu$ and $\Sigma_\mu(\nu)=+\infty$ otherwise.
By construction we have $\Sigma_V(\nu)\ppq \Sigma_{\mu_V}(\nu)$ with equality for all probability measures on $S_{\mu_V}$. 
Informally, an other way to express the link between the two is:
$$\Sigma_{\mu}=\sup_{V|\mu_V=\mu}\Sigma_V=\Sigma_{-2U_\mu+\infty{\bf\un}_{S_\mu^c}}.$$

With this new quantity, Theorem \ref{main} can be stated as follows:
\begin{theo}[Free  $T_1$ inequality, version for probability measures]
For any $\mu \in \p(\RR),$
  with  compact support such that its logarithmic potential
$ U_\mu(x) = -    \int \ln|x-y| d\mu(y)$ exists and is continuous on $\CC$,
there exists a constant $B_\mu$ such that
 for any probability measure $\nu$ 
$$ W_1^2(\nu, \mu) \ppq B_\mu \Sigma_\mu(\nu). $$
\end{theo}

But since $\mu$ is compactly supported the result of Bobkov and Götze also applies and gives:
$$ W_1^2(\nu, \mu) \ppq C_\mu H(\nu|\mu).$$

A natural question
is to ask whether our free inequality is a direct consequence of the classical one.
This is not the case thanks to the following:

\begin{proposition} 
Let $\lambda$ be the uniform law on $[0;1]$,
then
$$ \sup_{\nu \in \p([0,1]),\nu\neq \lambda} \frac{H(\nu|\lambda)}{\Sigma_\lambda(\nu)} =\infty.
$$
\end{proposition}

\begin{dem}
The proof of the property is essentially a direct calculation.
Consider $\nu_n$ the uniform law on
$$\bigcup_{i=0}^{n-1}\frac{i}{n}+\left[0;\frac{1}{n^2}\right].$$
Then $H(\nu_n)=\ln(n)$ but $\Sigma_\mu(\nu_n)$ remains bounded since the double logarithmic part is equivalent to the convergent Riemann sum
$$\frac{1}{n^2}\sum_{1\ppq i\neq j\ppq n}\ln\left(\frac{i}{n}-\frac{j}{n}\right).$$
\end{dem}

\subsection{Concentration property for $\beta$-ensembles}

As mentioned in our quick review of classical transport-entropy inequalities at the beginning of the introduction and detailed in
\cite{GL}, those inequalities are intimately linked with concentration properties of the measures involved.
Bolley, Guillin and Villani show in \cite{BGV} how to deduce from Talagrand's inequalities explicit bounds on the convergence
of the empirical measure of independent variables towards their common measure.
For example if $X_1,\dots,X_n,\dots$
are independent variables in $\RR^d$ of law $\mu$ satisfying $T_p(C)$ with $1\ppq p\ppq 2$, then for any $d'<d$, any $C'<C$,
there exists $N_0>0$ such that for all $N>N_0$, for all $\theta> v(N/N_0)^{-\frac{1}{2+d'}}$
$$\PP\left(W_1\left(\frac{1}{N}\sum_{i=1}^N\delta_{X_i},\mu\right)>\theta\right)<e^{-\gamma_p \frac{C'}{2}N\theta^2}$$
with $\gamma_p$ an explicit constant depending on $p$ in a very simple way. These results have been extended in \cite{Bo} and \cite{BoT}. 

Similarly, in our context, as we know that, under $\PNV,$ the empirical measure $\mun$ converges almost surely to $\mu_V,$
it is natural to
ask whether we can control the tail of the distribution of the random variable $W_1(\mun,\mu_V)$.

More generally,  we will deduce from Theorem \ref{main} a concentration result for the so-called $\beta$-ensembles for $\beta>0,$
that is for the empirical measure  of the $x_i$'s distributed according to the measure
$$\PNVb(dx_1,\dots,dx_N)= \prod_{i<j}|x_i -x_j|^\beta\exp\left(-N\sum_{i=1}^NV(x_i)\right)\frac{\prod_{i=1}^Ndx_i}{\ZNVb}.$$
This time the $x_i$ will asymptotically distribute according to the probability measure $\mu_{\frac{2V}{\beta}}.$

In comparison with Theorem \ref{main}, we need here some additional assumptions,
for technical reasons that will appear more clearly along the proofs.
Let us define $\|f\|_\Lip^A$  the Lipschitz norm of $f$ on a compact set $A$:
$$\|f\|_\Lip^A=\sup_{s,t\in A,s\neq t}\left|\frac{f(t)-f(s)}{t-s}\right|$$

\begin{hyp}\label{hyp2}\mbox{}\\
\vspace{-0.5cm}
\begin{enumerate}
 \item[a.] $V$ satisfies Hypothesis \ref{hyp}, is  locally Lipschitz, differentiable outside a compact set 
and there exists  $\alpha>0,d\pgq 2$ such that, 
$|V^\prime(x)|  \sim_{|x|\to+\infty} \alpha |x|^{d-1}.$
\item[b.] $V$ and $\beta >0$ are such that the equilibrium measure $\muVb$ has finite classical entropy. 
\end{enumerate}

\end{hyp}

The condition b. is not as restrictive as it may seem due to a result by Deift, Kriecherbauer and McLaughlin.
A direct consequence of the main result in \cite{DKM} is that this is satisfied as soon as $V$ is $\mathcal C^2.$ 
Note also that consequently 
Hypothesis \ref{hyp2} is satisfied in the particular case of a polynomial of even degree with positive leading 
coefficient.\par\medskip

Our concentration  result around the limiting measure is 
as follows~:

\begin{theo}[Concentration for $\beta$-ensembles] \label{theo:conc}
Let $V$ and $\beta >0$ satisfy Hypothesis \ref{hyp2}. Then there exists $u,v >0$  such that for any $\theta>v\sqrt{\frac{\ln (1+N)}{ N}}$,
$$ \PNVb\left(W_1(\mun,\mu_{\frac{2V}{\beta}} )\pgq
\theta\right) \ppq e^{-u N^2 \theta^2}.$$
\end{theo}

The result above is stated for potentials
which are  equivalent to a power at infinity but this hypothesis can be relaxed
if we restrict it to a compact set, as will be stated in Theorem \ref{conc:cpct}. \par\medskip

We want to emphasise that very few results are known in this direction.
Nevertheless, for matrix models ($\beta=1$ or 2) with strictly convex potentials, Proposition 4.4.26 in \cite{AGZ} shows 
that if $V$ is $\mathcal C^\infty$ with $V''\pgq \kappa>0$ and $V^\prime$ has a polynomial growth at infinity, 
then for all $\theta\pgq 0$, for all $1$-Lipschitz function $f$,
$$\PNVb\left(\left|\frac{1}{N}\tr f-\int\frac{1}{N}\tr fd\PNVb\right|>\theta\right)<e^{-N^2\frac{\kappa\theta^2}{2}}.$$

The strength of our result is that it is valid for any $\beta>0$, does not require any convexity assumption 
and gives a bound simultaneously on all Lipschitz functions
since $W_1(\mu,\nu)=\sup_{f\,\, 1-\Lip}|\mu(f)-\nu(f)|$.
On the other hand, our method does not allow to get a bound for all $\theta\pgq 0$ and the constant in the exponential decay is not explicit.\par\medskip

The rest of the paper is divided in two parts, the first one proves the free transport-entropy inequality Theorem \ref{main}; 
the second one deduces from there the concentration estimate Theorem \ref{theo:conc}.\par\medskip

\section{Free $T_1$  inequality}
This section is devoted to the proof of our main result Theorem \ref{main}.
In the first part, we will build some useful tools from potential theory.
Then, in the second part of this section,  we prove the result restricted  to a fixed compact.
The third part of this section  extends the result on measures whose support is arbitrary.

\subsection{Lipschitz perturbations of the potential} \label{sec:potential}
The first ingredient of the proof is to evaluate the distance between the equilibrium measures corresponding to two potentials
obtained from one another by a Lipschitz perturbation. Propositions \ref{perturbf} and \ref{confinement} 
will be particularly useful in the case when the perturbation 
is Lipschitz but we state them in a slightly more general context.

Before giving the statements of these propositions, we first need the following lemma, that uses crucially the properties of the Hilbert transform.
This should be classical but we did not find a proper reference and we give its proof for the sake of completeness.
We denote by $L^2(\RR)$ the set on functions such that $\int f^2(x) dx < \infty. $

\begin{lemma} \label{lem:hilbert}
Let $\mu$ be a compactly supported probability measure on $\RR$ whose logarithmic potential $U_{\mu}(x)=-\int\ln|x-y|d\mu(y)$ is continuous on $\CC$.
Then if $g$ is a continuously differentiable function on $\RR$ with compact support,
 $$\int_\RR g(x)d\mu(x)=\int_\RR (Hg)^\prime(x)U_\mu(x)dx$$
 where $H$ is the Hilbert transform: for $f \in L^2(\RR),$ $\forall x \in \RR,$
$$(Hf)(x)=\dashint \frac{f(y)}{x-y}dy := \lim_{\eps \downarrow 0} \int_{\RR \setminus [x-\eps, x+\eps]} \frac{f(y)}{x-y}dy.$$
\end{lemma}

The proof of the Lemma uses the following properties of the Hilbert transform, that can be found e.g. in the work of Riesz \cite{R} 
(in particular in paragraph 20)~:
\begin{property}\mbox{ }
\vspace{-0.2cm}
\label{riesz}
 \begin{enumerate}
  \item[a.] $H$ is an isometry on $L^2(\RR)$, $H^2=-id.$
\item[b.]
If $\psi$ is analytic in a neighbourhood of $\RR$, $H\Ima \psi=-\Rea \psi$ 
\item[c.] If $f$ is in $L^2(\RR),$  differentiable and such that $f^\prime$ is in $L^2(\RR),$ then $Hf$ is differentiable
and $ (Hf)^\prime =H(f^\prime).$ Moreover, if $f$ is continuously differentiable, so is $Hf.$
 \end{enumerate}

\end{property}

We now prove Lemma \ref{lem:hilbert}.
                         
\begin{dem} 
For $y>0$ and $g$ continuously differentiable on $\RR$ with compact support, we define
$$\phi(y):=\Ima\int_\RR  g(x) \int_\RR \frac{1}{\pi(x+iy-t)}d\mu(t) dx.$$
On one hand, 
if $X$ is of law $\mu$ and $\Gamma$ is an independent Cauchy variable we can rewrite $\phi$ as a convolution:
$$
\phi(y)=\iint g(x) \frac{y}{\pi((x-t)^2+y^2)}d\mu(t)dx=\E[g(X+y\Gamma)]$$

Therefore, by dominated convergence, $\phi(y)=\E[g(X)]+\eps(y),$
with $\eps(y)$ going to zero as $y$ goes to zero. Otherwise stated,
when $y$ goes to zero, $\phi(y)$ converges to $\int g(t)d\mu(t)$.

On the other hand, for any $y>0,$ the function $x \mapsto  \Ima \int \frac{1}{\pi(x+iy-t)}d\mu(t)$ is in $L^2(\RR)$
and, by Property \ref{riesz}.a. above, 
$$\phi(y)=\int (Hg)(x)H\left(\Ima\int \frac{1}{\pi(\cdot+iy-t)}d\mu(t)\right)(x)dx.$$
Thus, as $z \mapsto  \int \frac{1}{\pi(z-t)}d\mu(t)$ is analytic in a neighbourhood of $\RR$, by Property \ref{riesz}.b.,
$$\phi(y)=-\Rea\int (Hg)(x)\left(\int \frac{1}{\pi(x+iy-t)}d\mu(t)\right)dx.$$
Then, as $U_\mu$ is supposed to be continuous and $g$ continuously differentiable,
an integration by parts gives
$$\phi(y)=\Rea\int (Hg)^\prime(x)U_\mu(x+iy)dx$$
As $g$ is compactly supported, one can easily check that there exists $K>0$ such that for $x$ large enough, 
$|(Hg)^\prime(x)|\ppq \frac{K}{x^2}.$ As $\mu$ is  compactly supported, for $x$ large enough and any $y>0,$
$|U_\mu(x+iy)| \ppq K \ln x.$ Therefore, by dominated convergence, $\phi(y)$ converges to 
$ \int (Hg)^\prime(x)U_\mu(x)dx$ as $y$ goes to zero.
\end{dem}

We can now state the first perturbative estimate.
\begin{proposition}[Dependancy of the equilibrium measure in the potential] \label{perturbf}
For any $L>0$, there exists a finite constant $K_L$ such that,
for any $V,W$ satisfying Hypothesis \ref{hyp}, if $\mu_V$ and $\mu_W$ are probability measures on  $[-L;L]$ then
$$ W_1(\mu_V, \mu_W) \ppq K_L\osc(V-W).$$
with $\osc(f)=\sup_\RR f - \inf _\RR f$.
\end{proposition}

\begin{dem}
Our main tool for this proof is the use of the logarithmic potentials of the measures involved. We have already seen
in Theorem \ref{theo:st} that they are closely related. 
Corollary I.4.2 in \cite{ST} gives us a valuable estimate
$$ \|U_{\mu_V} - U_{\mu_W}\|_\infty \ppq  \|V-W\|_\infty.$$

We will also crucially use  a dual formulation for the distance $W_1.$ Indeed,
the Kantorovich-Rubinstein theorem (see e.g. Theorem 1.14 in \cite{Vil}) gives that
$$W_1(\mu_V, \mu_W)=\sup_g \mu_V(g)-\mu_{W}(g)$$
where the supremum is taken over the set of 1-Lipschitz function on $\RR.$
 
Note that the quantity $\mu_V(g)-\mu_{W}(g)$ does not change if we add a constant to $g$ or if we change the values of $g$ outside  $[-L;L]$.
This observation and a density argument show that
$$W_1(\mu_V, \mu_W)=\sup_{g\in\G}\mu_V(g)-\mu_{W}(g)$$
with $\G$ the set of $\mathcal C^1,$ compactly supported, 1-Lipschitz function on $\RR$, vanishing  outside of $[-2L;2L].$ 
Let $g$ be in $\G$,
according to Lemma \ref{lem:hilbert},
$$\mu_V(g)-\mu_{W}(g) =\int (H g)^\prime(x) (U_{\mu_V}(x) - U_{\mu_W}(x))dx.$$
Indeed, all the assumptions of Lemma \ref{lem:hilbert} are fulfilled, as we know from Theorem I.4.8 of \cite{ST}
that $U_{\mu_V}$ and $U_{\mu_W}$ are continuous on $\CC$ as soon as $V$ and $W$ are.
Now we cut this integral into two. On one hand, as $g \in \G,$ $\|g\|_\infty\ppq 2L,$ we have
\begin{align*}
&\left|\int_{|x|>2L+1} (Hg)^\prime(x) (U_{\mu_V} - U_{\mu_W})(x)\right|\\
&\ppq \|V-W\|_{\infty}\|g\|_\infty\int_{|x|>2L+1,|y|<2L}\frac{dxdy}{|x-y|^2}
\ppq K_L^1 \|V-W\|_{\infty}
\end{align*}
with $K_L^1=2L\int_{|x|>2L+1}\int_{|y|<2L}|x-y|^{-2}.$

On the other hand, by Cauchy-Schwarz inequality and using Properties \ref{riesz}.a. and c. of the Hilbert transform
\begin{align*}
\left|\int_{|x|<2L+1} (Hg)^\prime(x) (U_{\mu_V} - U_{\mu_W})(x)\right|&\ppq \|V-W\|_{\infty}(4L+2)^{1/2}\| H(g')\|_2\\
&=(4L+2)^{1/2}\| g'\|_2\|V-W\|_{\infty}\\
&\ppq (4L+2)\|V-W\|_{\infty}
\end{align*}

Finally, it is easy to check  that $\mu_V$ depends on $V$ only up to an additive constant, thus we can always translate 
$V$ such that $\|V-W\|_{\infty}=2\osc(V-W)$.
Thus we have proved
$$\mu_V(g)-\mu_{W}(g) \ppq K_L\osc(V-W)$$
with $K_L=2(K_L^1+4L+2).$
As $K_L$ does not depend on $g \in \G,$ taking the supremum for $g$ in $\G$ gives the result.

\end{dem}

The next step  is to show that given a Lipschitz function
$f$ on a given interval $[-L;L]$ we can extend the function outside this interval while keeping a control on the support of $\mu_{V+f}$ 
independently of $f$.
This property is rather technical but crucial since we will need to consider functions $f$  of arbitrary Lipschitz constant and a priori 
there is no way to control uniformly the support of
$\mu_{V+f}$.

\begin{proposition}[Confinement Lemma]
\label{confinement}
Let $V$ be a function satisfying Hypothesis \ref{hyp}.
For any $L>0,$ there exists $\tL>L$ depending only on $L$ and the potential $V$ such that for any  $u\pgq 0$, for any $u$-Lipschitz function $f$
on $[-L,L]$ one can find a function $\hf$  such that
\begin{enumerate}
\item $\hf$ is a bounded $u$-Lipschitz function on $\RR$
\item for all $|x|<L$, $\hf(x)=f(x)$
\item the support of $\mu_{V+\hf}$ is included in $[-\tL,\tL]$
\item $\osc(\hf) \ppq 2u \tL.$
\end{enumerate}
\end{proposition}

\begin{dem}
Let $V$ and $L$ be fixed and let $f$ be a $u$-Lipschitz function defined on $[-L,L].$
Again,  since $\mu_V$ depends on $V$ only up to an additive constant, one can always assume that $f(0)=uL$ (so that $f$ and the function 
$\hf$ we are going to define  both stay positive).

Let $\tL>L$ be a constant to be defined later.
Let us define $\hf$ as the biggest $u$-Lipschitz function which extends $f$ and is constant on components of $[-\tL,\tL]^c.$
More explicitly, we have
$$\hf(x)=
\left\{
\begin{array}{ll}
f(x) & \mbox{si } |x|\ppq L\\
f(L)+u(x-L) & \mbox{if } L\ppq x\ppq\tL\\
f(L)+u(\tL-L) & \mbox{if } x\pgq\tL\\
f(-L)-u(L+x) & \mbox{if } -\tL\ppq x\ppq-L\\
f(-L)-u(L-\tL) & \mbox{if } x\ppq-\tL.
\end{array}
\right.
$$
Our goal will be to find a constant $\tL,$ independent of $f$ and $u,$ such that $\hf$ (which depends on $\tL$) fulfils
the requirements of Proposition \ref{confinement}.\par\medskip

Since $V+\hf$ satisfies Hypothesis \ref{hyp}, 
the equilibrium measure $\mu_{V+\hf}$ is well defined.
Let us have a look at the value of the minimiser of the entropy functional $J_{V+\hf}$
(as defined in the introduction):  if 
$\lambda$ denotes the Lebesgue measure on $[0;1],$
$$c_{V+\hf} =\inf_{\nu \in \p(\RR)}J_{V+\hf}(\nu)\ppq J_{V+\hf}(\lambda).$$
Besides,
$$
J_{V+\hf}(\lambda)=\lambda(V)+\lambda(\hf)- \iint_{[0;1]^2}  \ln|x-y| dxdy\ppq\max_{[0;1]}V+(L+1)u- \frac{3}{4}
$$
Thus
$$c_{V+\hf} \ppq M_L(1+u) $$
with $M_L$ a constant only depending on $L$ and the potential $V.$

This estimate will allow us to find a bound on the support $S_{\mu_{V+\hf}}$ of $\mu_{V+\hf}$.
Indeed, define  $b=\sup\left\{|x|\in\RR / x \in S_{\mu_{V+\hf}}\right\}.$ 
Now we prove that a good choice of $\tL$ (depending only on $V$) such that $b>\tL$ leads to  contradiction.

Let us first assume that $\tL$ is chosen and $b>\tL$. From Theorem \ref{theo:st}, for any
$x$ in the support of $\mu_{V+\hf}$:
$$ V(x) + \hf(x) =  -2U_{V+\hf}(x) +C_{V+\hf}$$
and replacing $C_{V+\hf}$ by its expression given in  Theorem \ref{theo:st},
$$ V(x) + \hf(x) +\mu_{V+\hf}(V+\hf)=  -2U_{V+\hf}(x) + 2c_{V+\hf}$$
For any $x \in \RR,$ $-U_{V+\hf}(x) \ppq \ln(|x|+b)$.
On the other side, according to Hypothesis \ref{hyp}, there exists $\alpha >0$ and $\beta \in \RR$ (depending on $V$) such that 
for any $x \in \RR,$ $V(x) \pgq \alpha x^2 + \beta.$
Besides, as $f(0) = uL,$ $\hf(x)\pgq u|\min(|x|,\tL)-L|$, thus  $\mu_{V+\hf}(V+\hf)\pgq \beta$.
Putting these facts together, one gets for $|x|>\tL$ in the support of $\mu_{V+\hf}$:
$$ \alpha x^2 + \beta+u|\min(|x|,\tL)-L|+\beta\ppq 2\ln(|x|+b)+2M_L (1+u).$$
Now take a sequence of points in the support of $\mu_{V+\hf}$ converging in absolute value to $b.$ We get at the limit that:
$$ \alpha b^2 + 2\beta+ u(\tL-L)\ppq 2\ln(2 b)+2M_L (1+u).$$
There exists some $\gamma_V>1$ such that the function $\alpha x^2 - 2\ln(2 x)+2\beta-2M_L$ is strictly positive for $|x|>\gamma_V$.
Now choose $\tL>\gamma_V$,
since $b>\tL>1$, we get:
$$ u(\tL-L)<(\alpha b^2 +2\beta - 2\ln(2 b)-2M_L)+ u(\tL-L) \ppq 2 M_L u.$$

Then if we also choose $\tL>L+2M_L$ large enough, this leads a contradiction.
To sum up,  for this choice of $\tL$, 
we have proven that it is absurd to suppose that the support of $\mu_{V+\hf}$ is not in $[-\tL;\tL].$
Otherwise stated, $\hf$ satisfies the third point of the proposition.
The other points are trivially satisfied by construction.
\end{dem}

\subsection{Derivation of the theorem for measures on a given compact}\label{sec:bg}
The next step is to show a weak version of our main theorem, in the sense that the constant
in the inequality between Wasserstein distance and free entropy depends on the support of the measures 
under consideration.

\begin{proposition}[Free $T_1$ inequality on a compact]
 \label{theosup}
Let $V$ be a function satisfying Hypothesis \ref{hyp}.
For all $L>0,$ there exists a constant $B_{V,L},$ depending only on $L$ and $V,$
such that, for any probability measure $\nu$ with support in $[-L,L],$
$$ W_1(\nu, \mu_V)^2 \ppq B_{V,L} \Sigma_V(\nu).$$
\end{proposition}

\begin{dem}We can assume without loss of generality that $L$ is large enough for the support of $\mu_V$ to be inside $[-L;L]$.
We are going to use a duality argument. 
We first recall that
$$W_1(\nu, \mu_V) =  \sup_{\phi\,\, 1-Lip} \mu_V(\phi) - \nu(\phi).$$
We will first show Proposition \ref{theosup} for $\nu$ of the form $\mu_W,$ with support in $[-L,L]$ 
and with $W$ continuous and which coincides with $V$ outside a large compact set,
so that it satisfies Hypothesis \ref{hyp}.

Let $f$ be a $u$-Lipschitz function and $g= -\widetilde f,$ with $ \widetilde f$ defined as in Proposition \ref{confinement}. 
$$\mu_W(g)-\mu_V(g)-\Sigma_V(\mu_W)\ppq \sup_{\nu\in\p(\RR)}(\nu(g)-\mu_V(g)-\Sigma_V(\nu))$$
Note that since $g$ is equal to $-f$ on $[-L;L]$, the left hand side is just $\mu_V(f)-\mu_W(f)-\Sigma_V(\mu_W)$.

Let us control the right hand side: for any $\nu \in \p(\RR),$
$$\nu(g)-\mu_V(g)-\Sigma_V(\nu)=-J_{V-g}(\nu)+J_V(\mu_V)-\mu_V(g).$$
But since $J_{V-g}$ is minimal at $\mu_{V-g}$ and $J_V$ is minimal in $\mu_V$, for any $\nu \in \p(\RR),$ we have
\begin{align*}
 \nu(g)-\mu_V(g)-\Sigma_V(\nu)&\ppq J_V(\mu_{V-g}) -  J_{V-g}(\mu_{V-g}) -\mu_V(g)\\
 &=\mu_{V-g}(g)-\mu_V(g)\ppq |\widetilde f|_\Lip W_1(\mu_{V+\widetilde{f}},\mu_V).
 \end{align*}
By construction the support of $\mu_{V+\widetilde f}$ is inside $[-\tL;\tL]$, with $\tL$ as defined in Proposition \ref{confinement}, 
we can then apply  Proposition  \ref{perturbf}:
$$ \mu_V(f)-\mu_W(f)-\Sigma_V(\mu_W)\ppq | \widetilde f|_\Lip K_{\tL} osc(\widetilde f) \ppq u^2 \tL K_{\tL} .$$

Thus, there exists a constant $A_{V,L}$ such that for any  
 $\phi$ 1-Lipschitz and $u>0$, taking $f=u\phi,$ we have
$$ u(\mu_V(f)-\mu_W(f))-A_{V,L} u^2\ppq \Sigma_V(\mu_W).$$
We can take the supremum of this expression in $\phi$  and then in $u$  to get:
$$W_1(\mu_W,\mu_V)^2\ppq 4 A_{V,L} \Sigma_V(\mu_W).$$

We now conclude the proof by extending it to general probability measures $\mu$ supported on $[-L,L]$ (not just those of the form 
$\mu_W$ for  $W$ continuous and which coincides with $V$ outside a compact set).                                                                                                                                                                                                                                                                                                                                                                                                                                                                                                                                                                                                                                                                                                                                                                                                                                                                                                                                                                                                         
The idea, following \cite{HP} p.216, is that for any $\mu$ compactly supported there exists a sequence
of potential $W_\eps$ such that $\mu_{W_\eps}$ is sufficiently good approximation of $\mu$.
Indeed, we approximate $\mu$ by the sequence of probability measures $\mu*\lambda_\eps$
where $\lambda_\eps$ is the uniform law on $[0;\eps]$ and $*$ is the usual convolution of measures.
As the function $-2U_{\mu*\lambda_\eps}(x)$ grows logarithmically with $x$ whereas $V(x)$ grows at least quadratically,
one can choose $R_\eps > 2L$ such that if $|x| >R_\eps,$ $V(x) > -2U_{\mu*\lambda_\eps}(x).$

We now define
$$W_\eps(x)= \left\{
\begin{array}{ll}
-2U_{\mu*\lambda_\eps}(x), & \textrm{ if } |x| \ppq R_\eps, \\
\left(2- \frac{|x|}{R_\eps}\right) (-2U_{\mu*\lambda_\eps}(x)) + \left(\frac{|x|}{R_\eps}-1  \right)V(x),  & \textrm{ if }
R_\eps  \ppq |x| \ppq 2R_\eps, \\
V(x),  & \textrm{ if } |x| \pgq 2R_\eps.
\end{array}
\right.
$$
The result of the convolution is sufficiently smooth so that $U_{\mu*\lambda_\eps}$ is well defined and
continuous on $\CC$.
Thus for all $\epsilon>0$, $W_\eps$ is continuous and coincides with $V$ outside $[-2R_\eps, 2R_\eps].$
Besides, using  Theorem \ref{theo:st} we see that $\mu_{W_\eps}=\mu*\lambda_\eps$.
Since $\mu*\lambda_\eps$ is the equilibrium measures of the potential $W_\eps,$ 
we can apply the first part of the proof: for all $\eps>0$,
$$W_1(\mu*\lambda_\eps,\mu_V)^2\ppq 4 A_{V,L} \Sigma_V(\mu*\lambda_\eps).$$
Using the convexity, it is easy to check that $\lim_{\eps \to 0} \Sigma_V(\mu*\lambda_\eps)\ppq \Sigma_V(\mu)$.
As we know that $W_1$ is lower semi-continuous, we now take the limit when $\epsilon$ goes to $0$:
$$W_1(\mu,\mu_V)^2\ppq 4 A_{V,L} \Sigma_V(\mu).$$
\end{dem}

\subsection{Extension to non-compactly supported measure}\label{sec:confin}

To deduce Theorem \ref{main} from Proposition \ref{theosup}, we have to control what happens far from the support of $\mu_V.$
The idea is that, since $V$ grows faster than some $a x^2,$ if the support of $\mu$ is far from the support of $\mu_V,$ $\Sigma_V(\mu),$ 
which is growing like  $V$ should be much larger than $ W_1(\mu, \mu_V)^2 $ which is growing rather like $x^2.$ Therefore, it is enough to control
what happens in a vicinity of the support of $\mu_V$ and this case was treated in Proposition \ref{theosup}.

More precisely we have the following,
\begin{lemma} \label{scrountch}
Let $V$ be a function satisfying Hypothesis \ref{hyp}.
 There exists $\gamma_V>0$ and $R_V$ depending only on $V$ such that for any $\mu \in \p(\RR),$ there exists
$\widetilde \mu$ supported in $[-R_V,R_V]$ such that 
$$ \Sigma_V(\widetilde\mu) \ppq\Sigma_V(\mu)$$
 $$\gamma_V W_1(\mu, \widetilde\mu)^2 \ppq\Sigma_V(\mu).$$ 
\end{lemma}

We postpone the proof of the lemma to the end of this section and we first check that we can now get our main result (Theorem \ref{main}).

\begin{dem}
Let $\mu \in \p(\RR)$ and $\widetilde \mu$ corresponding to $\mu$ as in Lemma \ref{scrountch}. Then, using the triangular inequality 
and Proposition \ref{theosup}
\begin{align*}
 W_1(\mu, \mu_V)^2 & \ppq  2 W_1(\widetilde \mu, \mu_V)^2 + 2 W_1(\widetilde \mu,  \mu)^2 \\
& \ppq  2 B_{V,R_V} \Sigma_V(\widetilde \mu) + \frac{2}{\gamma_V}\Sigma_V(\mu)\\
&\ppq  2 \left(B_{V,R_V}+  \frac{1}{\gamma_V}\right)\Sigma_V(\mu). 
\end{align*}

\end{dem}

Finally, we prove Lemma \ref{scrountch}.

\begin{dem}
Let $R_V$ be a constant to be chosen later.
There exists $\alpha \in [0,1]$ such that $\mu = (1-\alpha) \mu_1 + \alpha \mu_2,$
with $\mu_1 \in \p([-R_V,R_V])$ and $\mu_2 \in \p([-R_V,R_V]^c).$
Then our definition for $\widetilde \mu$ is:
$$\widetilde \mu=  (1-\alpha) \mu_1 + \alpha \lambda$$
with $\lambda$ the Lebesgue measure\footnote{This precise choice of the Lebesgue measure here is not important, we just need a 
compactly supported measure of finite free entropy} on $[0;1].$ 

We now want to show the following statement, which implies both inequalities stated in the Lemma~: there exists $R_V$
 and $\gamma_V$ such that  
$$\Sigma_V(\mu) - \Sigma_V(\widetilde \mu) -\gamma_V W_1(\mu, \widetilde \mu)^2 \pgq 0.$$

Let us first bound the Wasserstein distance. In order to transport $\mu$ onto $\widetilde \mu$ one can always choose to 
transport $\mu_2$ to $\lambda$, this may not be optimal but gives the bound:
$$W_1(\mu, \widetilde \mu)^2\ppq\left(\alpha \int (|x|+1) d\mu_2(x)\right)^2 \ppq\alpha^2 \mu_2((1+|\cdot|)^2)$$

We then bound  the difference between  entropies:
\begin{align*}
 \Sigma_V(\mu) - \Sigma_V(\widetilde \mu)&\pgq
 \alpha (\mu_2-\lambda)(V) \\
 &-  \alpha^2 \iint \ln|x-y|[d\mu_2(x)d\mu_2(y) - d\lambda(x)d\lambda(y)]\\
&- 2\alpha(1-\alpha)  \iint \ln|x-y| d\mu_1(x)d(\mu_2-\lambda)(y) 
\end{align*}

Now we can get rid of the two double integrals by using that for all $x,y$, $\ln|x-y|\ppq\ln(1+|x|)+\ln(1+|y|)$ and 
that $|\int\ln|x-y|d\lambda(y)- \ln(1+|x|)|<C$ for some $C$ independent of $x$.
Thus,
\begin{align*}
\Sigma_V(\mu) - \Sigma_V(\widetilde \mu)&\pgq
\alpha( \mu_2(V- 2(1+\alpha)\ln(1+|\cdot|))\\
&-2\mu_1(\ln(1+|\cdot|))-\lambda(V-2(1+\alpha)\ln(1+|\cdot|))-2C).
\end{align*}

Finally, with $C_V=\lambda(V-4\ln(1+|\cdot|))+2C$, and the inequality $\mu_1(\ln(1+|\cdot|))\ppq \ln(1+R_V)\ppq \mu_2(\ln(1+|\cdot|))$,
$$ \Sigma_V(\mu) - \Sigma_V(\widetilde \mu)-\gamma_V W_1(\mu, \widetilde \mu)^2 \pgq
\alpha \mu_2(V- 6\ln(1+|\cdot|)-\alpha \gamma_V(1+|\cdot|)^2-C_V).
$$

We want this last expression to be positive.
We first choose $\gamma_V>0$ such that $\liminf_{|x| \rightarrow \infty} \frac{V(x)}{\gamma_V x^2}>1.$

Then $V(x)- 6\ln(1+|x|)-\gamma_V(1+|x|)^2-C_V$ goes to infinity when $|x|$ goes to infinity.
In particular we can choose $R_V>0$ such that for all $|x|>R_V$, it is positive.
Since $\mu_2$ has its support inside $[-R_V;R_V]^c$, the above expression is positive.
Since the choices of $\gamma_V$ and $R_V$ depend on $V$ only, this conclude the proof.
\end{dem}

\section{Concentration inequality for random matrices}
\label{sec:concentration}
In this section
we present an application of the free $T_1$ inequality to a result of concentration
for the empirical measure of a matrix model.
The concentration result holds not only on usual matrix models as defined in the introduction 
but also on the slightly more general family measures, usually called  $\beta$-ensembles.
We recall the definition of these models in the next section before proving our concentration estimates.

\subsection{$\beta$-ensembles}
For $\beta>0$, and $V$ a function satisfying Hypothesis \ref{hyp}, the $\beta$-ensemble with
potential $V$ is the family of laws on $\RR^N$, for $N>0,$ given by
$$\PNVb(dx_1,\dots,dx_N):= \prod_{i<j}|x_i -x_j|^\beta\exp\left(-N\sum_{i=1}^NV(x_i)\right)\frac{\prod_{i=1}^Ndx_i}{\ZNVb}.$$
with $\ZNVb$ a normalising constant which always exists under Hypothesis \ref{hyp}.
For $\beta=1,2,4$ this corresponds to the law of the eigenvalues of a matrix model (corresponding to the measure $\mu_V^N$
when $\beta=2$).

Some of the results stated in the introduction still holds for these models.
In particular, we  can still express nicely $\PNVb$ in terms of the empirical measure of the $x_i$'s. 
If $\mun:=\frac{1}{N}\sum_{i=1}^N\delta_{x_i}$ then
$$\PNVb(dx_1,\dots,dx_N)= \exp\left(-N^2\frac{\beta}{2}{\widetilde J_{\frac{2V}{\beta}}}(\mun)\right)\frac{\prod_{i=1}^Ndx_i}{\ZNVb}$$
with the functional $\widetilde J_V$ whose definition we recall
$$\widetilde J_V(\mu) = \int V(x)d\mu(x)-\iint_{x\neq y} \ln|x-y| d\mu(x) d\mu(y).$$
Similarly to the definition of $\Sigma_V,$ we also define
$$\widetilde \Sigma_V(\mu)= \widetilde J_V(\mu)- c_V.$$

One can expect that in the large $N$ limit, the eigenvalues  should organise this time according to the measure  $\mu_{\frac{2V}{\beta}}$.
This is indeed the case and we have a result analogous to Theorem \ref{theo:gd}, also proved in \cite{BaG}.
\begin{theo}[Large deviation principle for $\beta$-ensembles]\label{theoldpb}
Let $V$ be a function satisfying Hypothesis \ref{hyp}.
Under the law $\PNVb$
the sequence of random measures $\mun$ satisfies a large deviation principle in the speed $N^2$ with good rate function $\Sigma_{\frac{2V}{\beta}}$.
\end{theo}

\subsection{Approximate free $T_1$ inequalities for empirical measures}
At fixed $N$ we have to work with probability measure $\mun$ which have the drawback of being discrete.
This prevents us of applying the transport-entropy inequality since $\Sigma_V(\mun)=+\infty$.
We settle for an approximate inequality where  $\Sigma_V$ is replaced by $\StV$.

\begin{proposition}[Approximate free $T_1$ inequality]
\label{approxi}
Let $V$ be a locally Lipschitz function satisfying Hypothesis  \ref{hyp}. Then,
 for any $\mathcal{K}$ compact of $\RR,$ any $N\in \NN^*$ and any $(x_1,\dots,x_N) \in \mathcal{K}^N,$  
$$ W_1^2(\mun,\mu_V) \ppq 2B_V\StV(\mun)+3\frac{\|V\|^{\mathcal K_1}_\Lip+B+\ln(N)}{N}$$
where $B_V$ is the same constant as in Theorem \ref{main}, $B$  some universal finite constant and $\mathcal K_u$  the set 
of reals at distance less than $u$ from $\mathcal K.$
\end{proposition}

\begin{dem}
Let $\mathcal{K}$ be a compact set of $\RR$ and  $x_1,\dots,x_N$ be in $\mathcal K$.
The idea is to replace $\mun$ by a measure $\nun$ such that $W_1(\mun,\nun)$ is small
and $\SV(\nun)$ is close to $\StV(\mun)$.

We first  spread each $x_i$ such that they are at least $N^{-2}$ apart.
Let the $x_{(i)}$'s be the $x_i$'s rearranged by increasing order~: $x_{(1)}  \ppq x_{(2)}   \ppq \ldots \ppq x_{(n)},$
then define the $y_i$ by:
$$
\left\{
\begin{array}{lll}
y_1 &=& x_{(1)}\\
y_{i+1} &=& y_i + \max(x_{(i+1)}-x_{(i)}, \frac{1}{N^2})
\end{array}
\right.
$$

Then we define
$$\rhon=\frac{1}{N} \sum_{i=1}^N \delta_{y_i} \quad \textrm{ and } \quad \nun=\rhon*\lambda_{N^{-3}}$$
where $\lambda_{N^{-3}}$ is the uniform measure on $[0,N^{-3}]$ and $*$ is the usual convolution of measures.

Let us see how the Wasserstein distance and the entropy change when we replace $\mun$ by $\nun$.
Note that since $|y_i-x_{(i)}|<(i-1)N^{-2}$,
$$ W_1(\mun,  \nun)\ppq \frac{1}{N} \sum_{i=1}^N|y_i-x_{(i)}|\ppq \frac{1}{2N}$$
but
$$ W_1(\rhon,  \nun)\ppq\frac{1}{N^{3}}.$$
so that
$$W_1(\mun,  \nun)\ppq \frac{2}{N}.$$

Moreover, for any $i \neq j,$ $\ln |y_i-y_j| \pgq \ln|x_{(i)} - x_{(j)}|,$ and $y_i\in K_{N^{-1}}\subset K_1$,
$$ \StV(\mun)-\StV(\rhon) \pgq -\|V\|^{\mathcal K_1}_\Lip W_1(\mun,\rhon)\pgq -\|V\|^{\mathcal K_1}_\Lip \frac{2}{N}.$$

Let $(Z_i)_{i\pgq 1}$ and  $(\widetilde Z_i)_{i\pgq 1}$ be two independent families of independent variables uniformly distributed on $[0,1].$
We can express the difference of entropies using this variables:
\begin{multline*}
\StV(\rhon)-\SV(\nun) \pgq \int V(x) d(\rhon-\nun)(x) \\
+\frac{1}{N^2} \sum_{i \neq j} \E \left(\ln \left(1+ N^{-3} \frac{Z_i-Z_j}{y_i-y_j}\right)\right)
+\frac{1}{N^2} \sum_{i=1}^N \E\left(\ln N^{-3}|Z_i - \widetilde Z_i|\right)
\end{multline*}
Since for $i\neq j$, $|y_i-y_j|\pgq N^{-2}$, for $N>2$,
$$\E \left(\ln \left(1+ N^{-3} \frac{Z_i-Z_j}{y_i-y_j}\right)\right)\pgq \ln\left(1-\frac{2}{N}\right).$$
Thus,
$$\StV(\rhon)-\SV(\nun) \pgq -   \|V\|^{\mathcal K_1}_\Lip N^{-3}-\frac{B+3\ln N}{N}$$
with $B>0$ a finite constant.

This leads to, 
\begin{align*}
\StV(\mun)&=(\StV(\mun)-\StV(\rhon))+(\StV(\rhon)-\SV(\nun))+\SV(\nun)\\
&\pgq -3\frac{\|V\|^{\mathcal K_1}_\Lip+B+\ln(N)}{N}+\SV(\nun).
\end{align*}

Then by applying our the free transport inequality of Theorem \ref{main} for the potential $V$ on $\nun$, we obtain:
\begin{align*}
W_1(\mun,  \mu_V)^2&\ppq\left(W_1(\nun,  \mu_V)+\frac{2}{N}\right)^2\\
&\ppq 2W_1(\nun,  \mu_V)^2+\frac{8}{N^2}\\
&\ppq 2B_V\SV(\nun)+\frac{8}{N^2}\\
&\ppq 2B_V \StV(\mun)+\frac{8}{N^2}+3\frac{\|V\|^{\mathcal K_1}_\Lip+B+\ln(N)}{N}
\end{align*}

\end{dem}

\subsection{Tightness}
The next  step is to get a lower bound on the normalising constant $\ZNVb.$
From large deviation results (Theorem \ref{theoldpb}), it is easy to check
 that $\frac{1}{N^2} \ln \ZNVb$ has a finite limit $-\cVb$ and that
$\cVb=\frac{\beta}{2}c_{\frac{2V}{\beta}}$.
But hereafter, we are seeking a lower bound which is not asymptotic in $N.$ This is the only place where the condition b. of Hypothesis \ref{hyp2}
is needed.

\begin{lemma} \label{lbZ}
 For any $V$  a function satisfying Hypothesis \ref{hyp} and $\beta >0$ such that the equilibrium measure $\muVb$ 
such that $H(\muVb)$ is finite,
there exists a constant $\AVb$ such that for any $N \in \mathbb N^*,$
$$ \frac{1}{N^2} \ln \ZNVb  +\cVb\pgq\frac \AVb N.$$
\end{lemma}

\begin{dem}We follow closely a proof by Johansson in \cite{J}.\\
We denote by $\rho_V$ the density of $\mu_V.$
Note that if  $H(\muVb)$ is finite, it implies in particular that $\rVb$ is well defined
 and we introduce the following ensemble:
$$ E_N:= \left\{(x_1, \ldots, x_N) \in \RR^N | \prod_{i=1}^N \rVb(x_i) >0\right\}.$$

Then,
$$\ZNVb\pgq\int_{E_N}  
\exp\left(-N^2\frac{\beta}{2}{\widetilde J_{\frac{2V}{\beta}}}(\mun)\right)
\prod_{i=1}^Ne^{-\ln \rVb(x_i)}\prod_{i=1}^N\rVb(x_i)dx_i$$
and using Jensen inequality we get:
\begin{align*}
\ln \ZNVb&\pgq
-N^2\frac{\beta}{2}\int \widetilde J_{\frac{2V}{\beta}}(\mun)\prod_{i=1}^N\rVb(x_i)dx_i-N\int \ln \rVb(x)\rVb(x)dx\\
&=-N(N-1)\frac{\beta}{2}J_{\frac{2V}{\beta}}(\muVb)-N\int (V(x)+\ln \rVb(x))\rVb(x)dx
\end{align*}
We conclude by recalling that by definition: $J_{\frac{2V}{\beta}}(\muVb)=c_{\frac{2V}{\beta}}$.
\end{dem}

We then need to control the behaviour of the largest eigenvalue. The proof follows the ideas of the proof of Proposition 2.1 in \cite{BG2}. 
\begin{lemma}
 \label{tight}
Assume that $V$ is a continuous function such that
for some $\alpha >0$ and $d>1,$ $V(x)-\alpha x^d$ is bounded from below on $\RR$.
  Then for any $\beta>0$ and $0<a<\alpha/2$, there exists $M_0>0$ such that for any $M\pgq M_0$ and $N\in\NN^*$,
$$ \PNVb\left(\max_{i=1..N}|x_i| \pgq M\right) \ppq e^{-a NM^d}.$$
\end{lemma}
\begin{dem}
First, we need to control ${\ZNmVb}/{\ZNVb}$.
For all $L>0$,
\begin{multline*}
\frac{\ZNVb}{\ZNmVb}\pgq \int_{|x_N|<L} 
\int \exp\left({-(N-1)V(x_N)+\sum_{i=1}^{N-1} \ln |x_N-x_i|^\beta -V(x_i) }\right)\\
d\PNmVb(x_1, \ldots, x_{N-1})Y_{V,L}d\rho_{V,L}(x_N)
\end{multline*}
with $\rho_{V,L}$ the probability measure of density $(Y_{V,L})^{-1}\exp(-V(\cdot)){\bf \un}_{[-L;L]}$
and $Y_{V,L}$ its normalising constant.

By Jensen inequality we get
\begin{multline*}
 \ln\frac{\ZNVb}{\ZNmVb}\pgq \ln Y_{V,L} \\
+ \int\left(  -(N-1)V(x_N)+\sum_{i=1}^{N-1} \ln |x_N-x_i|^\beta -V(x_i)\right)\\
d\PNmVb(x_1, \ldots, x_{N-1})d\rho_{V,L}(x_N)
\end{multline*}

By Chebychev inequality, for any $R>0,$
$$\PNVb\left(\frac{1}{N}\sum_{i=1}^{N}V(x_i)>R\right)\ppq 
e^{- \frac 1 2 N^2R}\frac{Z^N_{V/2, \beta}}{\ZNVb}.$$
Now, from Theorem \ref{theoldpb}, we know that 
$ \frac{1}{N^2} \ln\left(\frac{Z^N_{V/2, \beta}}{\ZNVb}\right)$ converges so that it is bounded.
From there, we can easily deduce that $\int \frac{1}{N}\sum_{i=1}^{N}V(x_i)d\PNVb$
is uniformly bounded in $N.$
Since $x\mapsto \int \ln |y-x|d\rho_{V,L}(x)$ is bounded from below, we immediately see that there
exists a finite constant $\DVb$ such that for all $N$,
$$\frac{1}{N}\ln\frac{\ZNVb}{\ZNmVb}\pgq \DVb.$$ 

With this bound, we can complete the proof of  the Lemma.
We integrate separately on $x_N$ and on $x_1,\dots,x_{N-1}$ to get:
\begin{multline*}
\PNVb(|x_N|\pgq M)=\frac{\ZNmVb}{\ZNVb} \int_{|x_N|>M} e^{-NV(x_N)} 
\int \left(\prod_{i=1}^{N-1} |x_N-x_i|^\beta e^{-V(x_i)}\right)\\
d\PNmVb(x_1, \ldots, x_{N-1})dx_N.
\end{multline*}

There exists $b_{V,\beta}>0$ such that 
$$ |x-y|^\beta e^{-V(y)}\ppq b_{V,\beta}e^{V(x)/2}.$$
Therefore,
$$\PNVb(|x_N|\pgq M) \ppq
  e^{-ND_{V,\beta}}b_{V,\beta}^{N-1}\int_{|x_N|>M}e^{-\frac{N+1}{2}V(x_N)}dx_N,$$
Let $\gamma_V>0$ be such that for all $x$, $V(x)-\alpha x^d >-\gamma_V$. If $M>1$,
$$\PNV(|x_N|\pgq M) \ppq e^{-ND_{V,\beta}}b_{V,\beta}^{N-1}e^{\frac{N+1}{2}\gamma_V}2\frac{e^{-\frac{N+1}{2}\alpha M^d}}{\alpha\frac{N+1}{2}}.$$

For any $0<a< \alpha/2$, $M>M_0$, we obtain
$$ \PNVb(\max|x_i| \pgq M) \ppq N\PNV(|x_N|\pgq M)\ppq Ke^{-a NM^d}$$
with
$$K=\sup_{N\in\NN^*}Ne^{-ND_{V,\beta}}b_{V,\beta}^{N-1}e^{\frac{N+1}{2}\gamma_V}2\frac{e^{-(N(a-\frac{\alpha}{2})
+\frac{\alpha}{2}) M_0^d}}{\alpha\frac{N+1}{2}}.$$
Now, $a$ being fixed, we can clearly choose $M_0$ such that $K$ is finite and less than $1$.
\end{dem}

\subsection{Concentration results}
Our goal is now to show Theorem \ref{theo:conc}. As an intermediate result, we will first show the following result, which deals 
with concentration when restricted to a compact set. Then, the proof of Theorem \ref{theo:conc} will combine this result 
and the tightness shown in the preceding subsection.

\begin{theo}[Concentration inequality on a compact set]  
\label{conc:cpct}
Let $V$ be a locally Lipschitz function satisfying Hypothesis \ref{hyp} and $\beta>0,$ such that the equilibrium 
measure $\mu_{\frac{2V}{\beta}}$ has a finite classical entropy.
Then, for all $M>0$, there exists $u,v >0$  such that for all $\theta>v\sqrt{\frac{\ln (1+N)}{ N}}$,
$$ \PNVb\left(W_1(\mun,\mu_{\frac{2V}{\beta}} )\pgq
\theta,\forall i, |x_i|<M\right) \ppq e^{-u N^2 \theta^2}.$$
\end{theo}

\begin{dem}
We can rewrite our measure $\PNVb$ as follows
$$\PNVb(dx_1, \ldots, dx_N) = \frac{e^{-N^2 \cVb}}{\ZNVb}e^{- N^2 \frac{\beta}{2}\StVb(\mun)}dx_1 \ldots dx_N.$$

Thus, using Lemma \ref{lbZ}, we get
\begin{align*}&\PNVb\left(W_1(\mun, \muVb) \pgq \theta,\max|x_i|<M\right)\\
&\ppq e^{-N\AVb}(2M)^N \exp\left(-N^2\frac{\beta}{2}\inf\left\{\StVb(\mun)\left|\begin{array}{ll}\forall i, x_i\in[-M;M],\\ W_1(\mun, \muVb) \pgq \theta\end{array}\right.\right\}\right).
\end{align*}
Next we apply the approximate free $T_1$ inequality
of Proposition \ref{approxi} to obtain for any $u>0$,
\begin{multline*}
\PNVb\left(W_1(\mun, \muVb) \pgq \theta,\max|x_i|<M\right) \\
\ppq e^{-N\AVb}(2M)^N \exp\left(\frac{\beta N}{4B_V}\left(3(\|V\|^{[-M-1;M+1]}_\Lip+B+\ln(N))- N\theta^2\right)\right)\\
\ppq K(N,\theta,u)\exp\left(-u N^2 \theta^2\right)
\end{multline*}
 with
\begin{multline*}K(N,\theta,u) \\
= \exp\left(N\left(-\AVb+\ln(2M)+
\frac{3\beta}{4B_V}\left(\|V\|^{[-M-1;M+1]}_\Lip+B+\ln(N)\right)\right.\right.\\\left.\left.+
\left(u-\frac{\beta}{4B_V}\right)N\theta^2\right)\right).
\end{multline*}
Let us choose $u<\frac{\beta}{4B_V}$ so
 that $K(N,\theta,u)$ is a decreasing function in $\theta$.
 It is then easy to check that for a good choice of $v$ (which may depend on $M$, $V$ and $\beta$), for all $\theta>v\sqrt{\frac{\ln (1+N)}{ N}}$,
 $$K(N,\theta,u)\ppq K\left(N,v\sqrt{\frac{\ln (1+N)}{ N}},u\right)\ppq 1.$$ 
\end{dem}

We can now complete the proof of Theorem \ref{theo:conc}.

\begin{dem} 
Following the same steps as above, we get that for any $M, \theta >0,$
\begin{multline*}
 \PNVb\left(W_1(\mun, \muVb) \pgq \theta\right) \\
 \ppq e^{-N\AVb}(2M)^N \exp\left(\frac{\beta N}{4B_V}
\left(3(\|V\|^{[-M-1;M+1]}_\Lip+B+\ln(N))- N\theta^2\right)\right) \\
+\PNVb\left(\max|x_i|>M\right).
\end{multline*}

Now, from Lemma \ref{tight} above, under Hypothesis \ref{hyp2}, 
we have that, for any $0<a<\frac{\alpha}{2d}$ and $M$ large enough,
$$ \PNVb\left(\max|x_i|>M\right) \ppq e^{-aNM^{d}}.$$ 

Thus, if we  choose  $\theta>v\sqrt{\frac{\ln (1+N)}{ N}}$ with $v>M_0^{\frac{d}{2}}$ and $M=(\sqrt N\theta)^{\frac{2}{d}} > M_0$,
 we get, for any $u>0,$
$$\PNVb\left(W_1(\mun, \muVb) \pgq x\right)\ppq \widetilde K(N,\theta,u)\exp\left(-u N^2 \theta^2\right)$$
 with
\begin{multline*}
 \widetilde K(N,\theta,u) = 
\exp\left(N\left(-\AVb+\ln\left(2(\sqrt N \theta)^{\frac 2 d}\right)\right.\right.\\
+ \left.\left.
3\left(\|V\|^{[-(\sqrt N \theta)^{\frac{2}{d}}-1;(\sqrt N \theta)^{\frac{2}{d}}+1]}_\Lip+B+\ln(N)\right)+\left(u-\frac{\beta}{4B_V}\right)(\sqrt N\theta)^2\right)\right)\\
+\exp\left(-(a-u)N^2\theta^2\right).
\end{multline*}

 Again the result follows easily if we choose $u < \min\left(\frac{\beta}{4B_V} ,a\right)$
 since
 $$\|V\|^{[-(\sqrt N \theta)^{\frac{2}{d}}-1;(\sqrt N \theta)^{\frac{2}{d}}+1]}_\Lip=O((\sqrt N \theta)^{\frac{2(d-1)}{d}})=o(N\theta^2).$$
 
\end{dem}


\textbf{Aknowledgements:} We would like to thank Northeast Normal University in Changchun (China) for its hospitality
during the French-Chinese summer school in July 2011 where part of this work was completed. During this stay we could benefit 
from the help of Philippe Biane which allowed some substantial simplication of the proof of Proposition \ref{perturbf}.
We also thank François Bolley and Nathaël Gozlan for guiding us in the vast literature of optimal transport.\\

\bibliographystyle{alpha}
\bibliography{freeT1}

 \end{document}